# QUOTIENT CORRELATION: A SAMPLE BASED ALTERNATIVE TO PEARSON'S CORRELATION[1]

### By Zhengjun Zhang

### *University of Wisconsin*

The quotient correlation is defined here as an alternative to Pearson's correlation that is more intuitive and flexible in cases where the tail behavior of data is important. It measures nonlinear dependence where the regular correlation coefficient is generally not applicable. One of its most useful features is a test statistic that has high power when testing nonlinear dependence in cases where the Fisher's $Z$-transformation test may fail to reach a right conclusion. Unlike most asymptotic test statistics, which are either normal or $\chi^2$, this test statistic has a limiting gamma distribution (henceforth, *the gamma test statistic*). More than the common usages of correlation, the quotient correlation can easily and intuitively be adjusted to values at tails. This adjustment generates two new concepts—the tail quotient correlation and the tail independence test statistics, which are also gamma statistics. Due to the fact that there is no analogue of the correlation coefficient in extreme value theory, and there does not exist an efficient tail independence test statistic, these two new concepts may open up a new field of study. In addition, an alternative to Spearman's rank correlation, a rank based quotient correlation, is also defined. The advantages of using these new concepts are illustrated with simulated data and a real data analysis of internet traffic.

**1. Introduction.** In measuring the linear relation, Pearson's correlation, based on the product of $z$-scores, is extremely widely used. Applications related to correlation coefficients are thriving in many areas. However, there does not coexist a term with Pearson's correlation based on the quotient as far as we know, and furthermore, there is no analogue of the correlation

Received December 2006; revised April 2007.

[1]Supported in part by NSF Grants DMS-05-05528 and DMS-06-30210, Wisconsin Alumni Foundation, and Institute for Mathematical Research, ETH Zürich.

coefficient in extreme value theory [de Haan and de Ronde (1998)]. In this paper the quotient is chosen for its intuitive appeal.

This paper introduces a class of quotient correlations which can be used as an alternative to Pearson's correlation, and a class of rank based quotient correlations which can be used as an alternative to Spearman's rank correlation. Like both Pearson's correlation coefficients and Spearman's rank correlation coefficients, the quotient correlation coefficients and the rank based quotient correlation coefficients measure the level of agreement between two random variables.

In testing linear (nonlinear) independencies, although much work has been done, there is no public consensus on which method is the best. In the literature, Fisher (1921) introduces the extremely useful $Z$-transformation function, and Hotelling (1953) states that the best present-day usage in dealing with correlation coefficients is based on Fisher (1915, 1921). Recently, Freeman and Modarres (2005) showed that after Box–Cox transformation to approximate normal scales, Fisher's $Z$-transformation test is more efficient when testing independence between two positive random variables. Besides Fisher's $Z$-transformation test, many other testing procedures have been developed. These procedures are summarized in books by Cromwell, Labys and Terraza (1994), Cromwell, Hannan, Labys and Terraza (1994), Kanji (1999) and Thode (2002), among others.

We deliver a new test statistic, which we call the gamma test for independence, based on the quotient correlation concept. Considering that almost all existing asymptotic tests are normal or $\chi^2$ tests, our test is innovative because it has an asymptotic gamma distribution. More importantly, we show a theorem which assures that the limiting distribution of the test is still a gamma distribution even if one uses parametric estimated distribution functions or nonparametric empirical distribution functions to transform the data into unit Fréchet scales which are assumed in the definition of the test statistic, that is, parameter estimations and marginal transformations won't affect the limiting distribution of the test statistic. We shall show its efficiency and the inefficiency of Fisher's $Z$-transformation test when nonlinear dependence occurs.

There is a growing body of literature on modeling and testing tail dependence, also known as extremal dependence or asymptotic dependence, between the components of a two dimensional random vector, which refers to the concurrence of extreme values in each component. Examples dealing with asymptotically independent random variables can be found in Watson (1954), Sibuya (1960), de Haan and Resnick (1977), Chow and Teugels (1978), Anderson and Turkman (1991), Joe (1997) and Grady (2000). A significant step is due to Ledford and Tawn (1996, 1997). They introduced a sub-model which can be used to model tail dependence and nearly tail independence. Several papers, such as Peng (1999), Draisma, Drees, Ferreira and de Haan (2004), Ferreira and de Haan (2004) and others, have focused on estimation and rare events modeling, particularly of Dutch coast wind and water level extremes. However, it is not clear whether the data are tail dependent or not due to different conclusions. In the meantime, Hsing, Kluppelberg and Kuhn (2003) develop visualization methods to check tail dependence, and Falk and Michel (2004) study a different sub-model and graphic checking of tail dependence. These models are certainly useful and important in extreme value theory. However, constructing a test statistic for tail independence is still an interesting open problem [Campos, Marron, Resnick, Park and Jeffay (2003)]. Longin and Solnik (2001) and Campos et al. (2003) point out that in the study of extremal dependence, misleading results, even conflicting results, are often reported in the literature.

The quotient correlations can easily and intuitively be adjusted to values at tails. This adjustment enables us to construct tail dependence measures and to construct efficient test statistics for tail independence. This adjustment generates two new concepts—the tail quotient correlation and tail independence test statistics, which are also gamma statistics and effectively detect tail independencies or tail dependencies in all simulation examples to be presented here.

The paper is organized as follows: in Section 2 we introduce the quotient correlation coefficient and the rank based quotient correlation coefficient. We then in Section 3 develop the gamma test statistic for testing independence. We compare the gamma test with the Fisher's $Z$-transformation test. Section 4 introduces the tail quotient correlation concept. In Section 5 we obtain our gamma test statistics for testing tail independence. A necessary condition for two random variables being tail independent is presented here. In Section 6 simulation examples are analyzed. In Section 7 a real data example of internet traffic is tested. Concluding remarks are given in Section 8. Technically involved proofs are deferred to Section 9.

**2. The quotient.** Suppose $X$ and $Y$ are identically distributed positive random variables satisfying $\mathrm{P}(X \geq Y) > 0$, $\mathrm{P}(X \leq Y) > 0$. Then the quotients between $X$ and $Y$ are $Y/X$ and $X/Y$. Suppose $\{(X_i, Y_i),\ i = 1, \ldots, n\}$ is a bivariate random sample of $(X, Y)$. Then we will have $2n$ quotients.

2.1. *The quotient correlation coefficient.* It is easy to see that both $\max_{i \leq n}\{X_i/Y_i\} - 1$ and $\max_{i \leq n}\{Y_i/X_i\} - 1$ are asymptotically nonnegative. When these two values are both close to zero for a sufficiently large sample size $n$, one can say that the changing magnitudes of $X$ and $Y$ are very close to each other. The same is true for $\max_{i \leq n}\{X_i/Y_i\} \max_{i \leq n}\{Y_i/X_i\} - 1$ being close to zero.

We now use the above three quantities to define a new statistical coefficient, that is, the quotient correlation coefficient:

$$(2.1) \qquad \mathrm{q}_n = \frac{\max_{i \le n}\{Y_i/X_i\} + \max_{i \le n}\{X_i/Y_i\} - 2}{\max_{i \le n}\{Y_i/X_i\} \times \max_{i \le n}\{X_i/Y_i\} - 1}.$$

Some functional properties of $\mathrm{q}_n$ are illustrated by looking at the following function:

$$f(x, y) = \frac{x + y - 2}{xy - 1} \qquad \text{for } x \ge 1, \ y \ge 1, x + y > 2.$$

The case of $x + y = 2$ is trivial since it implies $X_i$ and $Y_i$ are equal. It is easy to see $0 \le f(x, y) \le 1$, $f(1, y) = 1$, $f(x, 1) = 1$, and $\frac{\partial f}{\partial x} = -\frac{(y-1)^2}{(xy-1)^2} < 0$, $\frac{\partial f}{\partial y} = -\frac{(x-1)^2}{(xy-1)^2} < 0$, $\frac{\partial^2 f}{\partial x \partial y} = -2\frac{(x-1)(y-1)}{(xy-1)^3} < 0$, which imply $f(x_1, \ y_1) \le f(x_2, \ y_2)$, where $x_1 \ge x_2$, $y_1 \ge y_2$.

The properties of the function $f(x, y)$ defined above suggest that, asymptotically, $\mathrm{q}_n$ can take values ranging from 0 to 1. The monotonicity of $f(x, y)$ suggests that, for a fixed sample size $n$, the larger the $\max_{i \le n}\{X_i/Y_i\}$ or the $\max_{i \le n}\{Y_i/X_i\}$, the smaller the $\mathrm{q}_n$, and hence the lower the agreement of changing magnitudes. The boundary properties of $f(x, y)$, that is, $f(1, y) = 1$, $f(x, 1) = 1$, suggest that the value of $\mathrm{q}_n$ is largely determined by the smaller value of $\max_{i \le n}\{X_i/Y_i\}$ and $\max_{i \le n}\{Y_i/X_i\}$. This property is very useful, especially when there exist outliers in one of the two sequences.

The computed $\mathrm{q}_n$ value can suggest the dependence magnitude between $X$ and $Y$. For example, when $\mathrm{q}_n$ is close to 1, there are co-movements between $X_i$ and $Y_i$, and hence, $X$ and $Y$ are nearly completely dependent.

Notice that in the definition of $\mathrm{q}_n$ the values of $X_i$ and $Y_i$ are required to be positive. In practice, this may not be the case. A variable transformation can resolve the issue. Sometimes a simple shift of the data to the right works too. For example, we can have $Y_i' = c - Y_i$, $i = 1, \ldots, n$, being positive values, where $c$ is chosen such that both $\max_{i \le n}\{Y_i'/X_i\}$ and $\max_{i \le n}\{X_i/Y_i'\}$ are greater than 1. Then we can use $Y_i'$ and $X_i$ to define $\mathrm{q}_n$ as we did in (2.1). When $\mathrm{q}_n$ is close to 1, the movement of $X_i$ is opposite to that of $Y_i$.

We now give an example that shows even if the supports of $X$ and $Y$ are $(0, \infty)$, $X/Y$ and $Y/X$ are bounded, and hence, $\mathrm{q}_n$ is greater than 0 regardless of sample size.

EXAMPLE 2.1. Suppose $Z_l$, $l = 1, \ldots, L$, are independent unit Fréchet random variables with distribution function $F(x) = \exp(-1/x)$, $x > 0$. Let $X = \max_{1 \le l \le L} \alpha_{l1} Z_l$, $Y = \max_{1 \le l \le L} \alpha_{l2} Z_l$, where $\alpha_{li} > 0$, $\sum_l \alpha_{li} = 1$, $i = 1, 2$. Let $c_1 = \max_{l \le L}(\alpha_{l1}/\alpha_{l2})$ and $c_2 = \max_{l \le L}(\alpha_{l2}/\alpha_{l1})$, then $\mathrm{P}(X/Y \le c_1) = \mathrm{P}(\max_l \alpha_{l1} Z_l < \max_l c_1 \alpha_{l2} Z_l) = 1$, and $\mathrm{P}(Y/X \le c_2) = \mathrm{P}(\max_l \alpha_{l2} Z_l < \max_l c_2 \alpha_{l1} Z_l) = 1$. Therefore, $\max_{i \le n}\{Y_i/X_i\} \xrightarrow{\text{a.s.}} c_1$, $\max_{i \le n}\{X_i/Y_i\} \xrightarrow{\text{a.s.}} c_2$, and $\mathrm{q}_n \xrightarrow{\text{a.s.}} \frac{c_1 + c_2 - 2}{c_1 * c_2 - 1}$, as $n \to \infty$.

This is a special case of multivariate maxima of moving maxima processes (or M4 processes) which have been relatively extensively studied in Smith and Weissman (1996), Zhang (2004, 2006), Zhang and Shinki (2006) and Zhang and Smith (2004a). This example will be used in Section 6.

Remark 2.2 stresses the importance of having the data used to compute $q_n$ in a unified scale.

REMARK 2.2. Suppose $\{(X_i^*, Y_i^*, U_i^*),\ i = 1, \ldots, n\}$ is a trivariate random sample of $(X^*, Y^*, U^*)$, where $X^*$, $Y^*$, $U^*$ are independent. Let $X_i = U_i^* X_i^*$, $Y_i = U_i^* Y_i^*$, then $X_i$ and $Y_i$ are dependent. However, $X_i/Y_i = X_i^*/Y_i^*$ which suggests that the value of $q_n$ computed based on $X_i$ and $Y_i$ and the value of $q_n$ computed based on $X_i^*$ and $Y_i^*$ are the same. At first glance the quotient correlation seemed unusable. Notice that $X_i$ and $X_i^*$ do not have the same distribution. When both $X_i$ and $X_i^*$ are converted into unit Fréchet scales, the quotients computed from the resulting variables will no longer give the same $q_n$ value. In other words, the quotient correlation coefficient is applicable.

It is natural to use unit Fréchet scores to define the quotient measure since the limiting distributions of maximum of quotients are a unit Fréchet distribution when each quotient is derived from two independent unit Fréchet (or unit exponential) random variables. When marginals are not unit Fréchet, marginal transformation is needed. Nonparametric methods of converting the data into unit Fréchet scales will be discussed in Section 2.2. For parametric method details, we refer to Smith (2003), Zhang (2005) and Zhang and Smith (2004b).

The quotient correlation $q_n$ measures the nonlinear dependence since the quotient itself is a nonlinear bivariate function. One can see that the value of $q_n$ mainly depends on the tail observations of two random variables $X$ and $Y$. If $X$ and $Y$ are tail independent, the value of $q_n$ turns out to be negligibly small. In this case, we propose to replace the maximal quotients by the $k$th largest quotients (or the $p$th percentiles) in their corresponding quotient series. In practice, $k$ should be small (1, 2, 3, etc.) in order to make the calculated dependence measures meaningful. Similarly, the choice of $p$ should be meaningful. These new statistics will be studied in a different project. We focus on (2.1) in this work.

2.2. *The rank based quotient correlation coefficient.* Suppose now that $X_{1n} \leq X_{2n} \leq \cdots \leq X_{nn}$ are the ordered sample of $X_1, X_2, \ldots, X_n$; $Y_{1n} \leq Y_{2n} \leq \cdots \leq Y_{nn}$ are the ordered sample of $Y_1, Y_2, \ldots, Y_n$. Define the empirical distribution function of the $X$ sample by $\widehat{F}_X(x) = 1/n \sum_{i=1}^n I_{(X_i \leq x)}$. Then $\widehat{F}_X(X_{kn}) = k/n$. Similarly, define the empirical distribution function of the $Y$ sample by $\widehat{F}_Y(y) = 1/n \sum_{i=1}^n I_{(Y_i \le y)}$. Then $\widehat{F}_Y(Y_{kn}) = k/n$. Therefore, when we transform the data into the unit Fréchet scale based on the empirical distribution functions and ranks, we get the same unit Fréchet scaled values for two identical ranks in both $X$ and $Y$ samples.

The transformation can be done with the function $-1/[\epsilon + \log\{\widehat{F}_X(X_{kn})\}] = -1/[\epsilon + \log(k/n)], k = 1, 2, \ldots, n$, where $\epsilon$ is a small number to avoid division by 0. In practice, however, this small number is sensitive to the transformed values depending on the sample size $n$. Similar discussions hold for the $Y$ sample too. Here, we introduce a different procedure to do data transformation.

Using the above arguments, we generate a unit Fréchet sample of size $n$ and denote $Z_{(1)} < Z_{(2)} < \cdots < Z_{(n)}$ as the ordered sample. The rank based quotient correlation is defined as

$$\mathrm{q}_n^{\mathfrak{R}} = \frac{\max_{i \le n}\{Z_{(\mathrm{rank}[Y_i])}/Z_{(\mathrm{rank}[X_i])}\} + \max_{i \le n}\{Z_{(\mathrm{rank}[X_i])}/Z_{(\mathrm{rank}[Y_i])}\} - 2}{\max_{i \le n}\{Z_{(\mathrm{rank}[Y_i])}/Z_{(\mathrm{rank}[X_i])}\} \times \max_{i \le n}\{Z_{(\mathrm{rank}[X_i])}/Z_{(\mathrm{rank}[Y_i])}\} - 1}.$$
(2.2)

One immediate advantage of this definition is that the transformation does not depend on any correction term. It is obvious that the definition is also based on a simulated random sample of a unit Fréchet random variable. We have found this form is easy to implement and the simulation results are very close to $\mathrm{q}_n$, which is based on a parametric transformation. In practice, we suggest simulating 10 or more unit Fréchet samples of size $n$, computing and reporting the mean or median value of $\mathrm{q}_n^{\mathfrak{R}}$.

After establishing quotient correlation coefficients, in the next section we present one of the important application features of quotient correlations, that is, their uses in testing independence between random variables.

## 3. The gamma test statistic for testing independence.

3.1. *The hypothesis and the test statistic.* It is well known that any absolutely continuous random variable can be converted into a random variable with desired properties by marginal transformation. Without loss of generality, suppose now $\{(X_i, Y_i), i = 1, \ldots, n\}$ is a bivariate random sample of $(X, Y)$ with unit Fréchet margins. Then we can formulate the following hypothesis test of independence:

$$H_0^c\colon X \text{ and } Y \text{ are independent} \quad \text{versus} \quad H_1^c\colon X \text{ and } Y \text{ are dependent.}$$
(3.1)

Here the superscript $c$ means that the complete data are used in (3.1). Later we shall give cases where only tail values are used.

The following theorem shows the properties of a random sample under $H_0^c$ of (3.1).

THEOREM 3.1. *If $X$ and $Y$ are independent and have unit Fréchet margins, and $(X_i,\ Y_i)$, $i = 1, \ldots, n$, is a random sample from $(X, Y)$, then $\lim_{n \to \infty} P\{n^{-1} \times (\max_{i \le n} Y_i/X_i + 1) \le x\} = e^{-1/x}$, and random variables $\max_{i \le n} Y_i/X_i$ and $\min_{i \le n} Y_i/X_i$ are asymptotically independent—that is,*

$$
\begin{aligned}
(3.2) \quad \lim_{n \to \infty} P\Big\{ &n^{-1}\Big( \max_{i \le n} Y_i/X_i + 1 \Big) \le x, \\
&n^{-1}\Big( \max_{i \le n} X_i/Y_i + 1 \Big) \le y \Big\} = e^{-1/x - 1/y}.
\end{aligned}
$$

*Furthermore, as $n \to \infty$, the random variable $q_n$ is asymptotically gamma distributed—that is,*

$$
(3.3) \qquad n q_n \xrightarrow{\mathcal{L}} \zeta,
$$

*where $\zeta$ is a* gamma$(2, 1)$ *random variable.*

A proof of Theorem 3.1 is deferred to Section 9.

The limiting distribution of $n q_n$ leads to the following testing procedure (the gamma test): if $n q_n > \zeta_\alpha$, where $\zeta_\alpha$ is the upper $\alpha$th percentile of the gamma$(2, 1)$ distribution, $H_0^c$ of (3.1) is rejected; otherwise, it is retained.

As mentioned earlier, we can transform the underlying random variables to any other scales of interest. This gives the power from a probabilistic point of view one may assume without loss of generality that the marginal distributions are known. Under the null hypothesis, although we do not have any parameters that associate one population with another, we do have parameters in each individual population. Estimators for the marginal distribution may influence the asymptotic behavior of tests. In a statistical problem $X_i$ could be a unit Fréchet random variable resulting from theoretical variable transformation, while $\tilde{X}_i(n)$ could be a unit Fréchet scaled value resulting from a variable transformation using an estimated distribution function. In most applications $\tilde{X}_i(n)$ are used in statistical inferences including testing tail independencies.

In the following Theorem 3.2, we introduce a sufficient condition under which the marginal effects from the estimated parameter values do not have any effects on the asymptotic behavior of tests.

THEOREM 3.2. *Suppose $\theta$, $\theta_i$, $i = 1, \ldots, n$, are absolutely continuous random variables with distribution function $G_\theta(x)$, and $\eta$, $\eta_i$, $i = 1, \ldots, n$, are absolutely continuous random variables with distribution function $G_\eta(x)$. Denote any estimated distribution function for $\theta$ by $\widehat{G}_\theta(x)$, and any estimated distribution function for $\eta$ by $\widehat{G}_\eta(x)$. Define $\tilde{X}_i = F^{-1}(\widehat{G}_\theta(\theta_i))$, $\tilde{Y}_i = F^{-1}(\widehat{G}_\eta(\eta_i))$, $i = 1, \ldots, n$, where $F(x)$ is the unit Fréchet distribution* *function. Suppose $\widehat{G}_\theta(\theta_i) \xrightarrow{\text{a.s.}} G_\theta(\theta_i)$, $\widehat{G}_\eta(\eta_i) \xrightarrow{\text{a.s.}} G_\eta(\eta_i)$, for all $i$, as $n \to \infty$. Then random variables $\max_{i \leq n} \tilde{Y}_i / \tilde{X}_i$ and $\max_{i \leq n} \tilde{X}_i / \tilde{Y}_i$ are each asymptotically gamma distributed, and they are asymptotically independent—that is,*

$$(3.4) \quad \begin{aligned} \lim_{n \to \infty} &\mathrm{P}\left\{ n^{-1}\left(\max_{i \leq n} \tilde{Y}_i / \tilde{X}_i + 1\right) \leq x, \ n^{-1}\left(\max_{i \leq n} \tilde{X}_i / \tilde{Y}_i + 1\right) \leq y \right\} \\ &= e^{-1/x - 1/y}. \end{aligned}$$

*Furthermore, as $n \to \infty$, the random variable,*

$$(3.5) \quad \tilde{\mathrm{q}}_n = \frac{\max_{i \leq n}\{\tilde{Y}_i / \tilde{X}_i\} + \max_{i \leq n}\{\tilde{X}_i / \tilde{Y}_i\} - 2}{\max_{i \leq n}\{\tilde{Y}_i / \tilde{X}_i\} \times \max_{i \leq n}\{\tilde{X}_i / \tilde{Y}_i\} - 1}$$

*is asymptotically gamma distributed—that is, $n\tilde{\mathrm{q}}_n \xrightarrow{\mathcal{L}} \zeta$, where $\zeta$ is a gamma$(2, 1)$ random variable.*

A proof of Theorem 3.2 is deferred to Section 9.

An immediate practical application is to apply the empirical distribution function in $\widehat{G}_\theta(\cdot)$ and $\widehat{G}_\eta(\cdot)$, and we have the following corollary.

COROLLARY 3.3. *Suppose $X_1, \ldots, X_n$ are observed variables with absolutely continuous distribution $F_X(x)$, and $Y_1, \ldots, Y_n$ are observed variables with absolutely continuous distribution $F_Y(y)$. Define*

$$\widehat{F}_X(x) = \frac{1}{n + 1/n} \sum_{i=1}^n I_{(X_i < x)}, \qquad \widehat{F}_Y(y) = \frac{1}{n + 1/n} \sum_{i=1}^n I_{(Y_i < y)}.$$

*Then, as $n \to \infty$ the empirical distribution transformation based quotient correlation*

$$(3.6) \quad \begin{aligned} \mathrm{q}_n^{\mathfrak{E}} = &\left( \max_{i \leq n}\{\log(\widehat{F}_X(X_i))/\log(\widehat{F}_Y(Y_i))\} \right. \\ &\quad + \max_{i \leq n}\{\log(\widehat{F}_Y(Y_i))/\log(\widehat{F}_X(X_i))\} - 2 \Bigg) \\ &\times \left( \max_{i \leq n}\{\log(\widehat{F}_X(X_i))/\log(\widehat{F}_Y(Y_i))\} \right. \\ &\quad \left. \times \max_{i \leq n}\{\log(\widehat{F}_Y(Y_i))/\log(\widehat{F}_X(X_i))\} - 1 \right)^{-1} \end{aligned}$$

*is asymptotically gamma distributed—that is, $n\mathrm{q}_n^{\mathfrak{E}} \xrightarrow{\mathcal{L}} \zeta$, where $\zeta$ is a gamma$(2, 1)$ random variable.*

Notice that $\mathrm{q}_n^{\mathfrak{E}}$ and $\mathrm{q}_n^{\mathfrak{R}}$ are asymptotically equal in distribution, and we also have $n\mathrm{q}_n^{\mathfrak{R}} \xrightarrow{\mathcal{L}} \zeta$, as $n \to \infty$. We apply $\mathrm{q}_n^{\mathfrak{R}}$ in this work.

3.2. *Power comparisons between the gamma test and Fisher's Z-transformation.* We now assume $X$ and $Y$ are normal random variables for the moment. Suppose $r_{XY}$ is the sample correlation coefficient of a bivariate random sample of size $n$. Then Fisher's Z-transformation is $W = \frac{1}{2}\log(\frac{1+r_{XY}}{1-r_{XY}})$, which satisfies $E(W) \cong \frac{1}{2}\log(\frac{1+\rho_{XY}}{1-\rho_{XY}})$, $\mathrm{Var}(W) \cong \frac{1}{n-3}$, and $Z = \frac{W-E(W)}{\sqrt{\mathrm{Var}(W)}} \xrightarrow{\mathcal{L}} N(0,1)$, which in turn allows us to test different hypotheses about the correlation—see details in papers by Fisher ([1915], [1921]), Hotelling ([1953]), Freeman and Modarres ([2005]) and a book by Härdle and Simar ([2003]). In our applications we test the null hypothesis of $\rho_{XY} = 0$ versus the alternative hypothesis of $\rho_{XY} \neq 0$.

When observed variables are actually independent, we have found that both Fisher's Z-transformation test and the gamma test have low Type I error rate (less than 5%), which, suggests that both tests are practically efficient. We now show the efficiency of the gamma test when variables are nonlinear dependent.

EXAMPLE 3.1. In this example, the correlation coefficient between two random variables is zero, while the two random variables are actually dependent. For example, suppose $X$ is standard normal, and let $Y = X^2$. It is easy to show $\mathrm{cov}(X,Y) = 0$. Like the previous example, we simulate 100 bivariate random samples of size $n$, $n = 25, \ldots, 100$. For each random sample, we conduct both Fisher's Z-transformation test and the gamma test with significance level $\alpha = 0.05$. At each sample size $n$, we compute the empirical power of each test. The performances are illustrated in the left panel of Figure [1]. Clearly, one can see that when sample size is larger than 37, the gamma test is superior, but Fisher's Z-transformation test is not applicable in this instance.

Needless to say, we are not proposing a test statistic to replace Fisher's Z-transformation test. In our point of view, these two tests can be used as two primary tests for testing independence since in most cases they will agree with each other and both tests are easy to implement. However, if the gamma test fails to suggest independence, while Fisher's Z-transformation test does suggest independence, the data may be actually drawn from two nonlinearly dependent random variables.

As a comparison, we performed the rank based quotient tests (replaced $q_n$ by $q_n^{\mathfrak{R}}$) to the example above. The results are plotted in the right panel in Figure [1]. One can see that its performance is just as good as (or slightly better than) the quotient correlation.

Notice that the power curve is not very smooth, and they may not be very accurate. For the purpose of demonstrating that the gamma test outperforms the Fisher's $Z$-transformation test in cases when nonlinear dependencies occur, we argue that Figure 1 may be sufficiently enough. To estimate the power curve much more accurate, one can use monotone smoothing which uses the strengths from neighboring points; see He and Shi (1998), for example, we will apply this method in a different study.

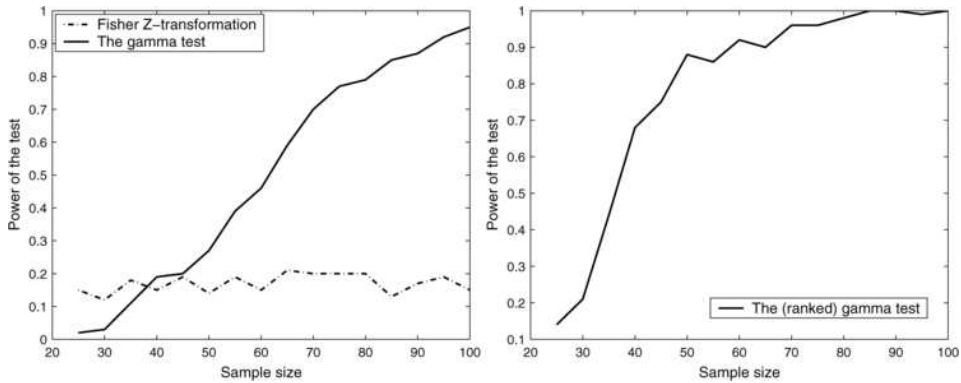

Fig. 1. *Comparison between Fisher's $Z$-transformation test and the gamma test. Left panel is based on $q_n$ and Fisher's $Z$-transformation function. Right panel is rank based test.*

**4. The tail quotient correlation.** In the literature the characterization of tail dependence between random variables is through a tail dependence index [Sibuya (1960) and Embrechts, McNeil and Straumann (2002)]. We present the definition next.

DEFINITION 4.1. Suppose random variables $X$ and $Y$ are identically distributed with $x_F = \sup\{x \in \mathbb{R} : P(X \le x) < 1\}$. The quantity (if exists): $\lambda = \lim_{u \to x_F} P(X > u | Y > u)$ is called the bivariate tail dependence index, or the tail dependence index, between $X$ and $Y$. It quantifies the amount of dependence of the bivariate upper tails of $X$ and $Y$. If $\lambda > 0$, then $X$ and $Y$ are tail dependent, otherwise the two random variables are tail independent.

This definition has been extended to cross-sectional tail dependencies in multivariate time series; see Zhang and Huang (2006) for example. Heffernan, Tawn and Zhang (2007) showed that $X$ and $Y$ are not necessarily identically distributed when one studies tail dependencies. However, we restrict our study to identically distributed random variables in this work.

To calculate $\lambda$, the joint distribution of $X$ and $Y$ must be known. In practice, when data is available, the data can be used to estimate the tail dependence index based on the empirical distribution functions and a very large threshold value $u$. We seek a new alternative quotient based statistical coefficient that measures the magnitude of tail dependence between two variables.

In general, when we study tail values, on one hand, we should use as much tail information as possible; on the other hand, we should use as little data from the center and below the threshold as possible. It is natural to use values above a threshold. However, this idea raises questions: how do we deal with the values below the threshold, and how do we get the quotient of the co-responses. Our choice is to use the threshold value for all those values below the threshold.

Using the arguments above, we define the tail quotient correlation coefficient as

$$q_{u,n} = \frac{\max_{i \le n}\{(u + W_i)/(u + V_i)\} + \max_{i \le n}\{(u + V_i)/(u + W_i)\} - 2}{\max_{i \le n}\{(u + W_i)/(u + V_i)\} \times \max_{i \le n}\{(u + V_i)/(u + W_i)\} - 1},$$
(4.1)

where $u$ is a positive threshold, and $W_i$ and $V_i$ are exceedance values over the threshold $u$ of positive random variables $X_i$ and $Y_i$ respectively. Here we assume $X_i$ and $Y_i$ are identically distributed.

Similarly, we can define a rank based tail quotient correlation coefficient. With the established notation in (2.2), a rank based tail quotient correlation is defined by

$$
\begin{aligned}
q_{u,n}^{\mathcal{R}} = {}& \left( \max_{i \le n}\left\{ \frac{u\mathbf{1}_{(W_i=0)} + Z_{(\operatorname{rank}[W_i])}\mathbf{1}_{(W_i>0)}}{u\mathbf{1}_{(V_i=0)} + Z_{(\operatorname{rank}[V_i])}\mathbf{1}_{(V_i>0)}} \right\} \right. \\
& + \max_{i \le n}\left\{ \frac{u\mathbf{1}_{(V_i=0)} + Z_{(\operatorname{rank}[V_i])}\mathbf{1}_{(V_i>0)}}{u\mathbf{1}_{(W_i=0)} + Z_{(\operatorname{rank}[W_i])}\mathbf{1}_{(W_i>0)}} \right\} - 2 \bigg) \\
& \times \left( \max_{i \le n}\left\{ \frac{u\mathbf{1}_{(W_i=0)} + Z_{(\operatorname{rank}[W_i])}\mathbf{1}_{(W_i>0)}}{u\mathbf{1}_{(V_i=0)} + Z_{(\operatorname{rank}[V_i])}\mathbf{1}_{(V_i>0)}} \right\} \right. \\
& \times \max_{i \le n}\left\{ \frac{u\mathbf{1}_{(V_i=0)} + Z_{(\operatorname{rank}[V_i])}\mathbf{1}_{(V_i>0)}}{u\mathbf{1}_{(W_i=0)} + Z_{(\operatorname{rank}[W_i])}\mathbf{1}_{(W_i>0)}} \right\} - 1 \bigg)^{-1},
\end{aligned}
$$
(4.2)

where $u$ is a global threshold value from a unit Fréchet distribution function and $\mathbf{1}_{(\cdot)}$ is an indicator function. We will use this rank based tail quotient correlation to conduct the tail independence test in Section 7.

One of the important uses of the tail quotient correlation is in the testing of tail independence which we study in the next section. Simulation examples are illustrated in Section 6.

**5. The gamma test statistics for testing tail independence.** The following theorem gives a necessary condition that two random variables are tail independent.

THEOREM 5.1 (Necessary condition). *In Definition* 4.1 *suppose that $X$ and $Y$ are unit Fréchet random variables, and $\lambda = 0$, then* $q_n \xrightarrow{a.s.} 0$, $q_{un} \xrightarrow{a.s.} 0$ *as $n \to \infty$.*

A proof of Theorem 5.1 is deferred to Section 9.

It is clear from Theorem 5.1 that a large tail quotient correlation would suggest tail dependence, and a zero value of tail quotient correlation would suggest that tail dependence is very weak. In most cases, however, the computed tail quotient correlation coefficient may never be zero. A statistical inference problem then arises: How small does a tail quotient correlation coefficient need to be for us to conclude that the two random variables are tail independent? This turns out to be a statistical hypothesis inference problem:

$$(5.1) \quad H_0^t \colon (X, Y) \text{ is tail independent} \quad \Leftrightarrow \quad H_1^t \colon (X, Y) \text{ is tail dependent}.$$

Here the superscript $t$ indicates that the test is for tail (in)dependence only.

Notice that (2.1) and (3.3) together test the null hypothesis of independence, not a null hypothesis of tail independence. We argue that the tail quotient correlation leads to a true tail independent test statistic under $H_0^t$ next.

Let

$$(5.2) \quad \begin{pmatrix} X_1, & X_2, & \ldots, & X_n \\ Y_1, & Y_2, & \ldots, & Y_n \end{pmatrix}$$

be an independent array of unit Fréchet random variables. Now let $(U_i, Q_i)$, $i = 1, \ldots, n$ be a bivariate random sequence, where both $U_i$ and $Q_i$ are correlated and have support over $(0, u]$ for a high threshold value $u$. Let $X_{ui} = X_i I_{\{X_i > u\}} + U_i I_{\{X_i \le u\}}$, $Y_{ui} = Y_i I_{\{Y_i > u\}} + Q_i I_{\{Y_i \le u\}}$, $i = 1, \ldots, n$. Then

$$(5.3) \quad \begin{pmatrix} X_{u1} \\ Y_{u1} \end{pmatrix}, \begin{pmatrix} X_{u2} \\ Y_{u2} \end{pmatrix}, \ldots, \begin{pmatrix} X_{un} \\ Y_{un} \end{pmatrix}$$

is a bivariate random sequence drawn from two dependent random variables $X_{ui}$ and $Y_{ui}$. Notice that $X_{ui} I_{\{X_{ui} > u\}}$ $(= X_i I_{\{X_i > u\}})$ and $Y_{ui} I_{\{Y_{ui} > u\}}$ $(= Y_i I_{\{Y_i > u\}})$ are independent, but $X_{ui} I_{\{X_{ui} \le u\}}$ $(= U_i I_{\{X_i \le u\}})$ and $Y_{ui} I_{\{Y_{ui} \le u\}}$ $(= Q_i I_{\{Y_i \le u\}})$ are dependent. In fact, one can easily construct arbitrarily dependent structure for values below the threshold $u$. Consequently, if only tail values are concerned, we can assume the tail values are drawn from (5.2) under the null hypothesis of tail independence. We have the following theorem.

THEOREM 5.2. *Suppose $V_i$ and $W_i$ are exceedance values (above $u$) in* (5.3), *and $U_i$ and $Q_i$ have the distribution $e^{1/u} e^{-1/x}$, $0 < x < u$. Then*

$$(5.4) \quad P\left( \frac{u + W_i}{u + V_i} \le t \right) = \begin{cases} \dfrac{t}{1+t} - \dfrac{t}{1+t} e^{-(1+t)/u}, & \text{if } 0 < t < 1, \\[2mm] \dfrac{t}{1+t} + \dfrac{1}{1+t} e^{-(1+t)/u}, & \text{if } t \ge 1, \end{cases}$$

$$(5.5) \quad \lim_{n \to \infty} P\left\{n^{-1}\left[\max_{i \le n}(u + W_i)/(u + V_i) + 1\right] \le x\right\} = e^{-(1 - e^{-1/u})/x}.$$

*Moreover,*

$$(5.6) \quad \lim_{n \to \infty} P\left\{n\left[\min_{i \le n}(u + W_i)/(u + V_i)\right] \le x\right\} = 1 - e^{-(1 - e^{-1/u})x}.$$

*The random variables* $\max_{i \le n}(u + W_i)/(u + V_i)$ *and* $\max_{i \le n}(u + V_i)/(u + W_i)$ *are tail independent—that is,*

$$(5.7) \quad \lim_{n \to \infty} P\left\{n^{-1}\left[\max_{i \le n}\frac{(u + W_i)}{(u + V_i)} + 1\right] \le x, \ n^{-1}\left[\max_{i \le n}\frac{(u + V_i)}{(u + W_i)} + 1\right] \le y\right\}$$
$$= e^{-(1 - e^{-1/u})/x - (1 - e^{-1/u})/y}.$$

*Furthermore, as* $n \to \infty$, *the random variable* $q_{u,n}$ *is asymptotically gamma distributed, that is,*

$$(5.8) \quad n q_{u,n} \xrightarrow{\mathcal{L}} \zeta,$$

*where* $\zeta$ *is* gamma$(2, 1 - e^{-1/u})$ *distributed.*

A proof of Theorem 5.2 is deferred to Section 9. Here are some intuitive interpretations. If $\max\{(u + W_i)/(u + V_i)\}$ or $\max\{(u + V_i)/(u + W_i)\}$ is close to 1, $q_{un}$ will have a larger value, which in turn implies that $n q_{un}$ has a larger probability to reject $H_0^t$ or a small $p$-value will be obtained. A larger $q_{un}$ value implies that there exists tail dependence (extreme co-movements) between $X_i$ and $Y_i$.

Equations (4.1) and (5.8) together establish a new test statistic—we call it the gamma test for the null hypothesis of tail independence—which can be used to determine whether there is tail dependence between two random variables. When $n q_{un} > \zeta_\alpha$, where $\zeta_\alpha$ is the upper $\alpha$th percentile of the gamma$(2, 1 - e^{-1/u})$ distribution, $H_0^t$ of (5.1) is rejected.

Considering the threshold value $u$ as a parameter, we obtain a new family (4.1) and (5.8) of gamma test statistics. When $u$ tends to zero, we obtain the special case of (2.1) and (3.3).

Notice that as $u$ is increasing, the value of $n q_{u,n}$ is increasing since both $\max\{(u + W_i)/(u + V_i)\}$ and $\max\{(u + V_i)/(u + W_i)\}$ are decreasing. Intuitively, a larger value of $n q_{u,n}$ renders to a rejection of the null hypothesis more likely. However, as $u$ is increasing, $1 - e^{-1/u}$ is decreasing, and a decreasing scale parameter $1 - e^{-1/u}$ leads to a larger cutoff value $\zeta_\alpha$. Therefore, a truncation (below the threshold value $u$) of data does not make a rejection of $H_0^t$ of (5.1) more likely. In fact, the gamma test can effectively detect all tail independent random variables in our simulation examples.

REMARK 5.1.    The tail quotient correlation coefficient $q_{un}$ can be used as a tail dependence measure once the null hypothesis of tail independence is rejected. In the literature several tail dependence measures have been proposed, such as those by Coles, Heffernan and Tawn (1999), Poon, Rockinger and Tawn (2001), Campos et al. (2003) and others. These measures and others have broad applications as long as tail events are concerned.

Also, notice that (5.8) holds when exceedance values of $V_i$ and $W_i$ are replaced by exceedance values of $\tilde{V}_i(n)$ and $\tilde{W}_i(n)$ resulting from $\tilde{X}_i(n)$ and $\tilde{Y}_i(n)$ respectively in Theorem 3.2.

The following theorem tells that one can apply the gamma tests when data are drawn from strict bivariate stationary processes which satisfy conditions specified in the theorem.

THEOREM 5.3.    *Suppose* $\{(X_i^*, Y_i^*), \ i = 1, \ldots, n\}$ *is a strict bivariate stationary processes with marginal distribution* $F_{X^*}(\cdot)$ *and* $F_{Y^*}(\cdot)$ *which are absolutely continuous.* $F(\cdot)$ *is a unit Fréchet distribution function. Let* $X_i = F^{-1}(F_{X^*}(X_i^*))$, $Y_i = F^{-1}(F_{Y^*}(Y_i^*))$, $i = 1, \ldots, n$. *Suppose the random sequences* $\{\xi_i = X_i/Y_i\}$, $\{\xi_i' = Y_i/X_i\}$ *satisfy Leadbetter's (1974)* $D(u_n)$ *and* $D'(u_n)$ *conditions for* $u_n = nx - 1$, *and Davis' (1982)* C1 *and* C2 *conditions. Suppose* $\widehat{F}_{X^*}(\cdot)$ *and* $\widehat{F}_{Y^*}(\cdot)$ *are estimated distribution functions for* $X_i^*$ *and* $Y_i^*$ *respectively, and as* $n \to \infty$, $\widehat{F}_{X^*}(X_i^*) \to F_{X^*}(X_i^*)$, $\widehat{F}_{Y^*}(Y_i^*) \to F_{Y^*}(Y_i^*)$, *for all* $i = 1, 2, \ldots$, *almost surely. Define* $\tilde{X}_i = F^{-1}(\widehat{F}_{X^*}(X_i^*))$ *and* $\tilde{Y}_i = F^{-1}(\widehat{F}_{Y^*}(Y_i^*))$. *Then as* $n \to \infty$, *both* (3.3) *and* (5.8) *hold, and* $n\tilde{q}_n \xrightarrow{\mathcal{L}} \zeta$, *where* $\tilde{q}_n$ *is defined as the one in* (3.5).

A proof of Theorem 5.3 is deferred to Section 9.

## 6. Simulation examples.

We have done substantial computation to check the efficiency of the gamma test. We have found that the gamma test has effectively detected tail (in)dependent random variables in all examples (including the sum and the maximum of a sequence of random variables, and variables in $M$-dependent sequences). Some typical examples are reported in the next section.

6.1. *Bivariate random samples.*    In this section we simulate bivariate random samples of sample size 500 from a bivariate random variable $(X, Y)$. The value of $\rho$ (correlation coefficient) is used to measure the correlation (linear dependence) between the two random variables. The value of $\theta$ in Gumbel type copulas is determined by $\theta = \sqrt{1-\rho}$.

Notice that the correlation coefficient between two random variables may not exist. If the correlation coefficient does not exist, $\rho$ does not have a specific meaning.

These simulated bivariate sequences are the following:

Table 1

*This table reports p-values obtained when the gamma test was applied to data drawn from models (a), (b), (c), (d), (e), (f), (g) and (h). The first column shows how the threshold values u were chosen*

| Percentile | (a) | (b) | (c) | (d) | (e) | (f) | (g) | (h) |
|---|---|---|---|---|---|---|---|---|
| 0.8000 | 0.2554 | 0.0000 | 0.0538 | 0.0000 | 0.6889 | 0.0538 | 0.0904 | 0.0000 |
| 0.8250 | 0.2448 | 0.0000 | 0.0512 | 0.0000 | 0.6860 | 0.0825 | 0.0906 | 0.0000 |
| 0.8500 | 0.2380 | 0.0000 | 0.0492 | 0.0000 | 0.6822 | 0.1226 | 0.0869 | 0.0000 |
| 0.8750 | 0.2340 | 0.0000 | 0.0473 | 0.0000 | 0.6769 | 0.1762 | 0.0834 | 0.0000 |
| 0.9000 | 0.2307 | 0.0000 | 0.0549 | 0.0000 | 0.6738 | 0.2228 | 0.1105 | 0.0000 |
| 0.9250 | 0.2277 | 0.0000 | 0.0745 | 0.0000 | 0.6668 | 0.2621 | 0.1962 | 0.0000 |
| 0.9500 | 0.2269 | 0.0001 | 0.0860 | 0.0000 | 0.6605 | 0.2624 | 0.1940 | 0.0000 |
| 0.9750 | 0.2308 | 0.0009 | 0.1010 | 0.0002 | 0.6596 | 0.2878 | 0.1975 | 0.0000 |

(a) Bivariate random samples drawn from two independent unit Fréchet random variables.

(b) Bivariate random samples drawn from the Gumbel copula,

$$C_\theta(H(u_1), H(u_2)) = e^{-(\sum_{i=1}^{2}[-\log\{H(u_i)\}]^{1/\theta})^\theta}, \qquad u_1 > 0, \ u_2 > 0,$$

where $H(u)$ is a unit Fréchet distribution function. Here $\theta = 0.4472$. Notice that the Gumbel copula is the only Archimedean copula which is also an extreme value copula [page 87, Mari and Kotz (2001)].

(c) Bivariate random samples drawn from the survival Gumbel copula,

$$P(X_{11} > u_1, X_{21} > u_2) = e^{-(\sum_{i=1}^{2}[-\log(P\{X_{i1} > u_i\})]^{1/\theta})^\theta},$$

$$u_1 > 0, \ u_2 > 0, \ \theta = 0.4472.$$

(d) Bivariate random samples drawn from Example 2.1, where $1 \le l \le 30$ and the coefficients are simulated values.

(e) Bivariate random samples drawn from $(1/U, \ 1/(1-U))$, where $U$ is uniform $(0, 1)$. (Sid Resnick contributed this example.)

(f) Bivariate normal random samples with correlation $\rho = 0.8$.

(g) Bivariate random samples drawn from $X = Z_1 E_1$, $Y = Z_2 E_1$, where $E_1$ is exponential 1. $Z_1$ and $Z_2$ are unit Fréchet.

(h) Bivariate random samples drawn from two t(4) (Student $t$ with 4 degrees of freedom) random variables with correlation $\rho = 0.8$.

For these eight examples, (a), (c), (e), (f) and (g) are tail independent examples, while (b), (d) and (h) are tail dependent examples. When we conduct the gamma test, the null hypotheses would be expected to be retained when data are drawn from any one of (a), (c), (e), (f) and (g), that is, larger $p$-values would be obtained; the null hypotheses would be expected to be rejected when data are drawn from any one of (b), (d) and (h), that is, smaller $p$-values would be obtained.

In our simulation study we do not use a fixed $u$ value. Instead, the threshold values are automatically chosen at the $100p$th percentiles of the sample data, where $p = 0.80, 0.825, 0.875, 0.90, 0.925, 0.95, 0.975$, respectively. We take the smallest value of the two $100p$th percentiles from two different sample sequences. Some typical test results from a single gamma test are summarized in Table 1. One can immediately see that the gamma test does work.

We have repeated running the program 100 times and recorded the total number of rejections of the null hypothesis of tail independence by the gamma test. The power of the gamma test among all simulated examples is 88%, which is very high. It would be fair enough to say that the gamma test is very efficient for testing tail independence.

**7. Application to internet traffic.** In network design, especially in internet traffic control, where the tail independence or dependence between two dependent variables can help networking researchers identify the main cause of internet delay, and design a better network structure, there have been a series of studies on tail dependence between internet traffic variables, for instance by Campos et al. (2003), Maulik, Resnick and Rootzén (2002), Resnick (2001, 2002), Zhang, Breslau, Paxson and Shenker (2002) and many others. These studies have shown that some variables are tail dependent, while some variables are tail independent.

Our application study tests tail (in)dependence of the joint behavior of large values of three internet traffic variables: size of response, time duration of response, and throughput (rate = size/time), which are in the context of HTTP (web browsing) responses. The data set was gathered from the University of North Carolina main link during April of 2001. The data set used in this analysis was for "morning" block (8:00AM-12:01PM). We obtained 76325 subsets of which there are 500 data points in each subset.

We performed the gamma test (using the rank based tail quotient correlation coefficient) on each subset. We found tail dependence between size of response and throughput in more than 3 percent of the data subsets. However, we did not find any statistical evidence that there is tail dependence between time duration and throughput. In Figure 2 we plot the tail behaviors between the traffic variables. The pattern in the right panel is discovered by the tail quotient test statistic. The data which generate the left panel is tested as tail independent. The data which generates the right panel is tested as tail dependent. As far as we know, such findings have not been reported in the literature.

These findings certainly encourage network researchers and network designers to better take into account the variable tail dependence in the network design and network traffic study. This is just a simple application of the gamma test in an internet traffic study. A detailed analysis of these traffic variables can be done, which will be a different research project. However, we restrict ourself to a traffic variable tail dependence test in this paper.

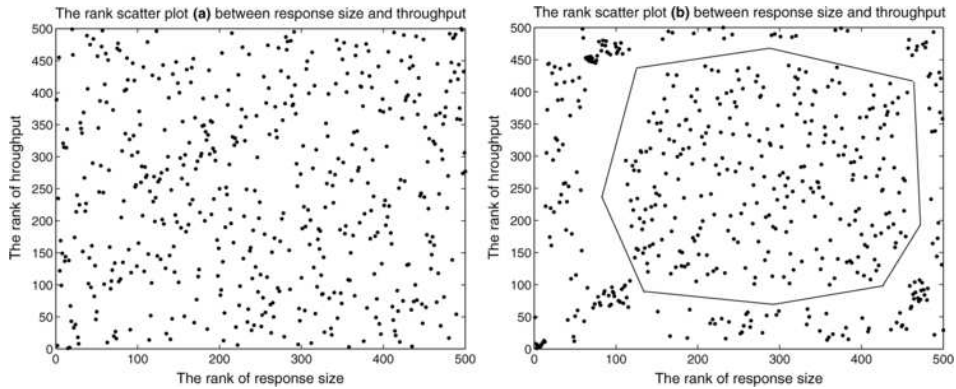

Fig. 2. *Comparison between two different data sets collected at different time periods. The pattern in the right panel is discovered by the tail quotient test statistic. The data which generates the left panel is tested as tail independent. The data which generates the right panel is tested as tail dependent.*

**8. Concluding remarks.** Resnick (1997) argues "Why non-linearities can ruin the heavy tailed modeler's day." What we hope to provide, in this paper, is an enabling methodology that will allow model builders to study tail behaviors of underlying multivariate time series on a scale previously not very attainable.

The new family of gamma test statistics can effectively lead to the right indications of tail (in)dependencies at high threshold levels, when the total sample size is about 500, as shown in our simulated examples. The gamma test has also given evidence of the tail independence of bivariate normal random variables, and of the tail independence of the sum and the maximum. The gamma test has outperformed Fisher's $Z$-transformation test when data are nonlinearly dependent and the sample size ranges from 37 to 100.

The new concept of the gamma test for (tail) independence should enhance the statistical theory of hypothesis testing. At the least, it adds one more member to the previously existing family of test statistics.

The tests and models of tail dependencies are still far from perfect. There is much work yet to be done. We believe that the methodology developed in this paper will be very useful in the analysis of financial data. Such research could produce more effective financial time series model and a clearer view of what has been missing in option pricing models.

In the literature, for exceedance problems, there have been existing test statistics derived based on empirical likelihood methods or empirical pro-cesses. Among those test statistics, the null hypothesis is tail dependencies which involve pre-specified model assumptions, and hence, they can be thought as testing a class of models. In our setup, the null hypothesis is tail independence. As a result, it is not practically useful to compare our test with existing tests, and it is not easy to carry out.

When the null hypothesis of (tail) independence is rejected, $q_n$ (or $q_{un}$) can be used as a (tail) dependence measure or as an estimate of a tail dependence index. The properties of this new measure may need to be studied further. A generalization of $q_{un}$ based on max-stable processes may be worth investigation.

The gamma tests and the (tail) dependence measures can be applied to many areas, for instance, financial study, internet traffic analysis (as mentioned above), environmental science, geophysics, microarray data analysis and psychological behavior study are just a few out of many areas.

## 9. Technical arguments.

PROOF OF THEOREM 3.1. First we have $P\{Y_i/X_i \le t\} = t/(1+t)$, which implies

$$\lim_{n \to \infty} P\left\{ n^{-1}\left( \max_{i \le n} Y_i/X_i + 1 \right) \le x \right\} = e^{-1/x}.$$

To prove (3.2), we compute the following probabilities:

$$P\left\{ n^{-1}\left( \max_{i \le n} Y_i/X_i + 1 \right) \le x, \ n^{-1}\left( \max_{i \le n} X_i/Y_i + 1 \right) \le y \right\}$$

$$= \left[ P\{X_i/(ny-1) < Y_i < (nx-1)X_i\} \right]^n$$

$$= \left[ \int_0^\infty \int_{z/(ny-1)}^{(nx-1)z} de^{-1/y} \, de^{-1/z} \right]^n = \left[ 1 - \frac{1}{nx} - \frac{1}{ny} \right]^n,$$

which implies (3.2). For (2.1), we have

$$nq_n = n \frac{\max\{Y_i/X_i\} + \max\{X_i/Y_i\} - 2}{\max\{Y_i/X_i\}\max\{X_i/Y_i\} - 1}$$

$$= \left( n^{-1}\left( \max\left\{ \frac{Y_i}{X_i} \right\} + 1 \right) + n^{-1}\left( \max\left\{ \frac{X_i}{Y_i} \right\} + 1 \right) - \frac{4}{n} \right)$$

$$\times \left( n^{-1}\left( \max\left\{ \frac{Y_i}{X_i} \right\} + 1 \right) n^{-1}\left( \max\left\{ \frac{X_i}{Y_i} \right\} + 1 \right) \right.$$

$$\left. - n^{-1}\left[ n^{-1}\left( \max\left\{ \frac{Y_i}{X_i} \right\} + 1 \right) + n^{-1}\left( \max\left\{ \frac{X_i}{Y_i} \right\} + 1 \right) \right] \right)^{-1}$$

$$= \frac{[n^{-1}(\max\{Y_i/X_i\} + 1) + n^{-1}(\max\{X_i/Y_i\} + 1)][1 + o_p(1)]}{[n^{-1}(\max\{Y_i/X_i\} + 1)n^{-1}(\max\{X_i/Y_i\} + 1)][1 + o_p(1)]},$$

which gives (2.1) by the Cramér–Wold device and Slutsky's theorem. $\square$

LEMMA 9.1. *Suppose $X$, $X_1, X_2, \ldots$, are positive random variables. Then $X_n \xrightarrow{\text{a.s.}} X$ if and only if there are two sequences of random variables $\xi_1(n)$, $\xi_2(n)$ such that $\xi_1(n) < X_n/X < \xi_2(n)$, $n = 1, 2, \ldots$, and $\xi_1(n) \xrightarrow{\text{a.s.}} 1$, $\xi_2(n) \xrightarrow{\text{a.s.}} 1$, as $n \to \infty$.*

PROOF. The sufficient condition is obvious. The necessary condition can be shown by simply defining $\xi_1(n) = X_n/X * (1 - 1/n)$, and $\xi_2(n) = X_n/X * (1 + 1/n)$. $\square$

PROOF OF THEOREM 3.2. Denote $X_i = F^{-1}(G_\theta(\theta_i))$, $Y_i = F^{-1}(G_\eta(\eta_i))$, then $\tilde{X}_i \xrightarrow{\text{a.s.}} X_i$, $\tilde{Y}_i \xrightarrow{\text{a.s.}} Y_i$, as $n \to \infty$ and, hence, by Lemma 9.1, there exist $\xi_j(n) > 0$, $\xi_j(n) \xrightarrow{\text{a.s.}} 1$, $j = 1, 2, 3, 4$, as $n \to \infty$, and

$$\xi_1(n) < \frac{\tilde{X}_i}{X_i} < \xi_2(n), \qquad \xi_3(n) < \frac{\tilde{Y}_i}{Y_i} < \xi_4(n), \qquad i = 1, \ldots, n,$$

which implies

$$P\left\{ \max_{i \le n} \frac{Y_i}{X_i} \frac{\xi_4(n)}{\xi_1(n)} \le nx - 1, \max_{i \le n} \frac{X_i}{Y_i} \frac{\xi_2(n)}{\xi_3(n)} \le ny - 1 \right\}$$

$$\le P\left\{ \max_{i \le n} \frac{\tilde{Y}_i}{\tilde{X}_i} \le nx - 1, \max_{i \le n} \frac{\tilde{X}_i}{\tilde{Y}_i} \le ny - 1 \right\}$$

$$\le P\left\{ \max_{i \le n} \frac{Y_i}{X_i} \frac{\xi_3(n)}{\xi_2(n)} \le nx - 1, \max_{i \le n} \frac{X_i}{Y_i} \frac{\xi_1(n)}{\xi_4(n)} \le ny - 1 \right\}.$$

Since $\xi_j(n)/\xi_k(n)$ converges to 1 almost surely for all $j$ and $k$, and $\max_{i \le n} \frac{Y_i}{X_i}$ and $\max_{i \le n} \frac{X_i}{Y_i}$ are asymptotically independent, so by Slutsky's theorems, we have both the first probability and the last probability in the above inequalities converging to $e^{-1/x - 1/y}$ as $n \to \infty$, and hence, (3.4) is true. The proof of the asymptotic distribution of $n\tilde{q}_n$ is similar to the proof of the asymptotic distribution $nq_n$ in Theorem 3.1. $\square$

PROOF OF THEOREM 5.1. We prove $q_n \xrightarrow{\text{a.s.}} 0$, as $n \to \infty$ only. The proof of $q_{un} \xrightarrow{\text{a.s.}} 0$ is similar.

Suppose without loss of generality that there is a finite number $z$, $1 < z < \infty$, such that $P(X/Y < z) = 1$. We have

$$\frac{P(X > u, Y > u)}{P(Y > u)} = \frac{P(zY > X > u, zY > zu)}{P(X > u)}$$

$$= \frac{P(zY > X > zu, zY > zu)}{P(X > u)}$$

$$+ \frac{\mathrm{P}(zY > X > u, zY > zu, u < X < zu)}{\mathrm{P}(X > u)}$$

$$\geq \frac{\mathrm{P}(X > zu)}{\mathrm{P}(X > u)} \to \frac{1}{z} > 0, \qquad \text{as } u \to \infty,$$

which contradicts with the zero limit on the left-hand side. So both ratios $X/Y$ and $Y/X$ have support over $(0, \infty)$. So $\max_{i \leq n} \frac{Y_i}{X_i} \xrightarrow{\text{a.s.}} \infty$, $\max_{i \leq n} \frac{X_i}{Y_i} \xrightarrow{\text{a.s.}} \infty$, as $n \to \infty$, and hence, $q_n \xrightarrow{\text{a.s.}} 0$. $\square$

PROOF OF THEOREM 5.2.   We first prove (5.4). For $t < 1$, we have

$$\mathrm{P}\left\{ \frac{u + W_i}{u + V_i} < t \right\}$$

$$= \mathrm{P}\left\{ \frac{u + W_1}{u + V_1} < t \right\} = \mathrm{P}\left\{ \frac{u + (Z_1 - u)I_{Z_1 > u}}{u + (Z_2 - u)I_{Z_2 > u}} < t \right\}$$

$$= \mathrm{P}\left\{ \frac{Z_1}{Z_2} < t, Z_1 > u, Z_2 > u \right\} + \mathrm{P}\left\{ \frac{u}{Z_2} < t, Z_1 < u, Z_2 > u \right\}$$

$$+ \mathrm{P}\left\{ \frac{Z_1}{u} < t, Z_1 > u, Z_2 < u \right\} + \mathrm{P}\left\{ 1 < t, \ Z_1 < u, \ Z_2 < u \right\}$$

$$= \triangle_1 + \triangle_2 + \triangle_3 + \triangle_4,$$

where $\triangle_3 = 0$, $\triangle_4 = 0$ and

$$\triangle_1 = \mathrm{P}\{u < Z_1 < tZ_2\}$$

$$= \int_u^\infty \int_{z_1/t}^\infty d e^{-1/z_2} \, d e^{-1/z_1}$$

$$= \int_u^\infty (1 - e^{-t/z_1}) \frac{1}{z_1^2} e^{-1/z_1} \, dz_1$$

$$= 1 - e^{-1/u} - \int_u^\infty \frac{1}{z_1^2} e^{-(1+t)/z_1} \, dz_1 = 1 - e^{-1/u} - \frac{1}{1+t}[1 - e^{-(1+t)/u}],$$

$$\triangle_2 = \mathrm{P}\{Z_1 < u, Z_2 > u/t\} = e^{-1/u}(1 - e^{-t/u}) = e^{-1/u} - e^{-(1+t)/u}.$$

Adding $\triangle_1$ and $\triangle_2$, we get (5.4) for $0 < t < 1$.

For $t \geq 1$, we have

$$\triangle_1 = \mathrm{P}\left\{ \frac{u}{t} < \frac{Z_1}{t} < Z_2, Z_2 > u \right\}$$

$$= \int_u^{tu} \int_u^\infty d e^{-1/z_2} d e^{-1/z_1} + \int_{tu}^\infty \int_{z_1/t}^\infty d e^{-1/z_2} \, d e^{-1/z_1}$$

$$= \frac{t}{1+t} + e^{-2/u} - e^{-1/u} - \frac{t}{1+t} e^{-(1+t)/(tu)},$$

$\triangle_2 = \mathrm{P}\{Z_1 < u, Z_2 > u\} = e^{-1/u}(1 - e^{-1/u}) = e^{-1/u} - e^{-2/u},$

$\triangle_3 = \mathrm{P}\{u < Z_1 < tu, Z_2 < u\} = e^{-1/u}(e^{-1/(tu)} - e^{-1/u}) = e^{-(1+t)/(tu)} - e^{-2/u},$

$\triangle_4 = \mathrm{P}\{Z_1 < u, \ Z_2 < u\} = e^{-2/u}.$

Adding the above terms together, we thus obtain (5.4).

We now prove (5.5). Since

$$\mathrm{P}\left\{\frac{u + W_1}{u + V_1} < nx - 1, \frac{u + V_1}{u + W_1} < ny - 1\right\}$$

$$= \mathrm{P}\left\{\frac{u + W_1}{u + V_1} < nx - 1, \frac{u + V_1}{u + W_1} < ny - 1, X_1 > u, Y_1 > u\right\}$$

$$+ \mathrm{P}\left\{\frac{u + W_1}{u + V_1} < nx - 1, \frac{u + V_1}{u + W_1} < ny - 1, X_1 > u, Y_1 < u\right\}$$

$$+ \mathrm{P}\left\{\frac{u + W_1}{u + V_1} < nx - 1, \frac{u + V_1}{u + W_1} < ny - 1, X_1 < u, Y_1 > u\right\}$$

$$+ \mathrm{P}\left\{\frac{u + W_1}{u + V_1} < nx - 1, \frac{u + V_1}{u + W_1} < ny - 1, X_1 < u, Y_1 < u\right\}$$

$$= \triangle_1 + \triangle_2 + \triangle_3 + \triangle_4,$$

where

$$\triangle_1 = \mathrm{P}\left\{\frac{Y_1}{X_1} < nx - 1, \frac{X_1}{Y_1} < ny - 1, X_1 > u, Y_1 > u\right\}$$

$$= \int_u^\infty \int_{\max\{x_1/(ny-1), u\}}^{(nx-1)x_1} de^{-1/y_1} de^{-1/x_1}$$

$$= \int_u^{(ny-1)u} \int_u^{(nx-1)x_1} de^{-1/y_1} de^{-1/x_1} + \int_{(ny-1)u}^\infty \int_{x_1/(ny-1)}^{(nx-1)x_1} de^{-1/y_1} de^{-1/x_1}$$

$$= \frac{nx-1}{nx} - \frac{nx-1}{nx} e^{-(nx)/(nx-1)u} + e^{-2/u} - \frac{1}{ny} - \frac{ny-1}{ny} e^{-(ny)/(ny-1)u},$$

$$\triangle_2 = \mathrm{P}\left\{\frac{u}{X_1} < nx - 1, \frac{X_1}{u} < ny - 1, X_1 > u, Y_1 < u\right\}$$

$$= \mathrm{P}\{u < X_1 < (ny-1)u, \ Y_1 < u\} = e^{-ny/(ny-1)u} - e^{-2/u},$$

$$\triangle_3 = \mathrm{P}\left\{\frac{Y_1}{u} < nx - 1, \ \frac{u}{Y_1} < ny - 1, \ X_1 < u, Y_1 > u\right\} = e^{-nx/(nx-1)u} - e^{-2/u},$$

and $\triangle_4 = e^{-2/u}$. Therefore,

$$\mathrm{P}\left\{\frac{u + W_1}{u + V_1} < nx - 1, \frac{u + V_1}{u + W_1} < ny - 1\right\}$$

$$= \triangle_1 + \triangle_2 + \triangle_3 + \triangle_4$$

$$= 1 - \frac{1}{nx} - \frac{1}{ny} + \frac{1}{nx} e^{-nx/(nx-1)u} + \frac{1}{ny} e^{-ny/(ny-1)u}.$$

Since

$$P\left\{ n^{-1} \left[ \max_{i \leq n} \frac{(u + W_i)}{(u + V_i)} + 1 \right] \leq x, \ n^{-1} \left[ \max_{i \leq n} \frac{(u + V_i)}{(u + W_i)} + 1 \right] \leq y \right\}$$

$$= P^n \left\{ \frac{u + W_1}{u + V_1} < nx - 1, \ \frac{u + V_1}{u + W_1} < ny - 1 \right\}$$

$$= \left[ 1 - \frac{1}{nx} - \frac{1}{ny} + \frac{1}{nx} e^{-nx/(nx-1)u} + \frac{1}{ny} e^{-ny/(ny-1)u} \right]^n$$

$$\to e^{-(1 - e^{-1/u})/x - (1 - e^{-1/u})/y} \qquad \text{as } n \to \infty,$$

which gives (5.5)–(5.7).

For (4.1) and (5.8), letting $\Upsilon_i = (u + W_i)/(u + V_i)$, $\Theta_i = (u + V_i)/(u + W_i)$, we have

$$n\mathsf{q}_{un} = n \frac{\max\{(u + W_i)/(u + V_i)\} + \max\{(u + V_i)/(u + W_i)\} - 2}{\max\{(u + W_i)/(u + V_i)\} \max\{(u + V_i)/(u + W_i)\} - 1}$$

$$= (n^{-1}(\max\{\Upsilon_i\} + 1) + n^{-1}(\max\{\Theta_i\} + 1) - 4/n)$$

$$\quad \times (n^{-1}(\max\{\Upsilon_i\} + 1) n^{-1}(\max\{\Theta_i\} + 1)$$

$$\quad \quad - n^{-1}[n^{-1}(\max\{\Upsilon_i\} + 1) + n^{-1}(\max\{\Theta_i\} + 1)])^{-1}$$

$$= \frac{[n^{-1}(\max\{\Upsilon_i\} + 1) + n^{-1}(\max\{\Theta_i\} + 1)][1 + o_p(1)]}{[n^{-1}(\max\{\Upsilon_i\} + 1) n^{-1}(\max\{\Theta_i\} + 1)][1 + o_p(1)]},$$

which gives (4.1) by the Cramér–Wold device and Slutsky's theorem. So the proof of the theorem is complete.  $\square$

PROOF OF THEOREM 5.3.  Using Theorem 3.5.2 of Leadbetter, Lindgren and Rootzén (1983) and Proposition 3.1 of Davis (1982), we immediately get $M_n = \max(\xi_1, \xi_2, \dots, \xi_n)$ and $W_n = \min(\xi_1, \xi_2, \dots, \xi_n)$ are asymptotically independent. Then the proof of (3.3) follows from Theorem 3.1, and the proof of (5.8) follows from Theorem 5.2. The last part follows from Theorem 3.2. $\square$

**Acknowledgments.**  The author thanks the referee, the Associate Editor, and the editor for their valuable constructive suggestions which improve the presentation of the paper.

## REFERENCES

ANDERSON, C. W. and TURKMAN, K. F. (1991). The joint limit distribution of sum and maxima of stationary sequences. *J. Appl. Probab.* **28** 33–44. MR1090445

CAMPOS, F. H., MARRON, J. S., RESNICK, S. I., PARK, C. and JEFFAY, K. (2003). Extremal dependence: Internet traffic applications. Available at http://www.cs.unc.edu/Research/dirt/proj/marron/ExtremalDependence/v2/LargeValAssoc2.pdf

CHOW, T. L. and TEUGELS, J. L. (1978). The sum and the maximum of i.i.d. random variables. In *Proc. Second Prague Symp. Asymptotic Statistics* 81–92. North-Holland, Amsterdam. MR0571177

COLES, S., HEFFERNAN, J. and TAWN, J. (1999). Dependence measures for extreme value analyses. *Extremes* **2** 339–365.

CROMWELL, J. B., LABYS, W. C. and TERRAZA, M. (1994). *Univariate Tests for Time Series Models*. Sage Publications, Thousand Oaks, CA.

CROMWELL, J. B., HANNAN, M. J., LABYS, W. C. and TERRAZA, M. (1994). *Multivariate Tests for Time Series Models*. Sage Publications, Thousand Oaks, CA.

DAVIS, R. A. (1982). Limit laws for the maximum and minimum of stationary sequences. *Z. Wahrsch. Verw. Gebiete* **61** 31–42. MR0671241

DE HAAN, L. and RESNICK, S. I. (1977). Limit theory for multivariate sample extremes. *Z. Wahrsch. Verw. Gebiete* **40** 317–337. MR0478290

DE HAAN, L. and DE RONDE, J. (1998). Sea and wind: Multivariate extremes at work. *Extremes* **1** 7–45. MR1652944

DRAISMA, G., DREES, H., FERREIRA, A. and DE HAAN, L. (2004). Bivariate tail estimation: Dependence in asymptotic independence. *Bernoulli* **10** 251–280. MR2046774

EMBRECHTS, P., MCNEIL, A. and STRAUMANN, D. (2002). Correlation and dependence in risk management: Properties and pitfalls. In *Risk Management*: *Value at Risk and Beyond* (M. A. H. Dempster, ed.) 176–223. Cambridge Univ. Press. MR1892190

FALK, M. and MICHEL, R. (2004). Testing for tail independence in multivariate extreme value models. In *Dependence Modeling*: *Statistical Theory and Applications in Finance and Insurance*. Quebec, Canada.

FERRERIA, A. and DE HAAN, L. (2004). On the consistency of the estimation of the probability of a failure set. In *The 3rd International Symposium on Extreme Value Analysis*: *Theory and Practice*. Aveiro, Portugal.

FISHER, R. A. (1915). Frequency distribution of the values of the correlation coefficient in samples from an indefinitely large population. *Biometrika* **10** 507–521.

FISHER, R. A. (1921). On the "probable error" of a coefficient of correlation deduced from a small sample. *Metron* **1** 1–32.

FREEMAN, J. and MODARRES, R. (2005). Efficiency of test for independence after Box–Cox transformation. *J. Multivariate Anal.* **95** 107–118. MR2164125

GRADY, A. (2000). A higher order expansion for the joint density of the sum and the maximum with applications to the estimation of the climatological trends. Ph.D. dissertation, Dept. Statistics, Univ. North Carolina.

HÄRDLE, W. and SIMAR, L. (2003). *Applied Multivariate Statistical Analysis*. Springer, Berlin. MR2061627

HE, X. and SHI, P. (1998). Monotone B-spline smoothing. *J. Amer. Statist. Assoc.* **93** 643–650. MR1631345

HEFFERNAN, J. E., TAWN, J. A. and ZHANG, Z. (2007). Asymptotically (in)dependent multivariate maxima of moving maxima processes. *Extremes* **10** 57–82.

HOTELLING, H. (1953). New light on the correlation coefficient and its transform. *J. Roy. Statist. Soc. Ser. B* **15** 193–232. MR0060794

HSING, T., KLÜPPELBERG, C. and KUHN, G. (2003). Modelling, estimation and visualization of multivariate dependence for risk management. Web reference.

JOE, H. (1997). *Multivariate Models and Dependence Concepts.* Chapman and Hall, London. MR1462613

KANJI, G. K. (1999). *100 Statistical Tests.* Sage Publications.

LEADBETTER, M. R. (1974). On extreme values in stationary sequences. *Z. Wahrsch. Verw. Gebiete* **28** 289–303. MR0362465

LEADBETTER, M. R., LINDGREN, G. and ROOTZÉN, H. (1983). *Extremes and Related Properties of Random Sequences and Processes.* Springer, Berlin. MR0691492

LEDFORD, A. W. and TAWN, J. A. (1996). Statistics for near independence in multivariate extreme values. *Biometrika* **83** 169–187. MR1399163

LEDFORD, A. W. and TAWN, J. A. (1997). Modelling dependence within joint tail regions. *J. Roy. Statist. Soc. Ser. B* **59** 475–499. MR1440592

LONGIN, F. and SOLNIK, B. (2001). Extreme correlation of international equity markets. *J. Finance* **LVI** 649–675.

MARI, D. D. and KOTZ, S. (2001). *Correlation and Dependence.* Imperial College Press.

MAULIK, K., RESNICK, S. and ROOTZEN, H. (2002). A network traffic model with random transmission rate. *J. Adv. App. Probab.* **39** 671–699. MR1938164

PENG, L. (1999). Estimation of the coefficient of tail dependence in bivariate extremes. *Statist. Probab. Lett.* **43** 399–409. MR1707950

POON, S.-H., ROCKINGER, M. and TAWN, J. (2001). New extreme-value dependence measures and finance applications. Univ. Strathclyde, 100 Cathedral Street, Glasgow G40LN, UK.

RESNICK, S. I. (1997). Why non-linearities can ruin the heavy tailed modeler's day. In *A Practical Guide to Heavy Tails*: *Statistical Techniques for Analyzing Heavy Tailed Distributions* (R. Adler, R. Feldman and M. Taqqu, eds). Birkhäuser, Boston. MR1652283

RESNICK, S. (2001). Modeling data networks. *SEMSTAT*, *Seminaire de Statistique 2001, Gothenberg, Sweden.* Chapman and Hall, London.

RESNICK, S. (2002). Hidden regular variation, second order regular variation and asymptotic variation, unpublished manuscript. Web reference.

SIBUYA, M. (1960). Bivariate extreme statistics. I. *Ann. Inst. Statist. Math.* **11** 195–210. MR0115241

SMITH, R. L. (2003). Statistics of extremes, with applications in the environment, insurance and finance. In *Extreme Value in Finance, Telecommunications and the Environment* (B. Finkenstadt and H. Rootzén, eds). Chapman and Hall/CRC, London.

SMITH, R. L. and WEISSMAN, I. (1996). Characterization and estimation of the multivariate extremal index. Manuscript, UNC.

THODE, H. C. JR. (2002). *Testing for Normality.* Dekker, New York. MR1989476

WATSON, G. S. (1954). Extreme values in samples from $m$-dependent stationary stochastic processes. *Ann. Math. Statist.* **25** 798–800. MR0065122

ZHANG, Y., BRESLAU, L., PAXSON, V. and SHENKER, S. (2002). On the characteristics and origins of internet flow rates. In *Proceedings of SIGCOMM02.* ACM, Pittsburgh.

ZHANG, Z. (2004). Some results on approximating max-stable processes. Submitted for publication.

ZHANG, Z. (2005). A new class of tail-dependent time series models and its applications in financial time series. *Advances in Econometrics* **20** 323–358.

ZHANG, Z. (2006). The estimation of M4 processes with geometric moving patterns. *Ann. Instit. Stat. Math.* 18 October 2006. Available at http://www.springerlink.com/content/21x436446679737m/.

ZHANG, Z. and HUANG, J. (2006). Extremal financial risk model and portfolio evaluation. *Comput. Statist. Data Anal.* **51** 2313–2338. MR2307503

ZHANG, Z. and SHINKI, K. (2006). Extreme co-movements and extreme impacts in high frequency data in finance. *J. Banking and Finance* **31** 1399–1415.

ZHANG, Z. and SMITH, R. L. (2004a). The behavior of multivariate maxima of moving maxima processes. *J. Appl. Probab.* **41** 1113–1123. MR2122805

ZHANG, Z. and SMITH, R. L. (2004b). On the estimation and application of Max-stable processes. Unpublished manuscript, Univ. Wisconsin.

DEPARTMENT OF STATISTICS
UNIVERSITY OF WISCONSIN
MADISON, WISCONSIN 53706
USA
E-MAIL: zjz@stat.wisc.edu